\documentclass{amsart}

\theoremstyle{plain}
\newtheorem{lemma}{Lemma}[section]

\newtheorem*{mainthm}{Main Theorem}

\newtheorem{con}[lemma]{Conjecture}

\theoremstyle{definition}
\newtheorem{defi}[lemma]{Definition}

\newtheorem{example}[lemma]{Example}

\theoremstyle{remark}
\newtheorem{remark}[lemma]{Remark}
\newtheorem{notation}[lemma]{Notation}

\newcommand{\p}{\mathbb{P}}

\newcommand{\z}{\mathbb{Z}}
\newcommand{\f}{\mathbb{F}}

\newcommand{\rea}{\mathbb{R}}
\newcommand{\oc}{{\mathcal O}}

\newcommand{\ls}{{\mathcal L}}

\newcommand{\I}{{\mathcal I}}
\newcommand{\e}{\epsilon}

\DeclareMathOperator{\h}{H}

\numberwithin{equation}{section}
\numberwithin{figure}{section}
\numberwithin{table}{section}

\begin{document}

\title{Special linear Systems on Toric Varieties}
\author{Antonio Laface and Luca Ugaglia}
\address{Dipartimento di Matematica, Universit\`a degli Studi di Milano, 
Via Saldini 50, 20133 MILANO}
\email{laface@mat.unimi.it\\
ugaglia@mat.unimi.it}
\begin{abstract}
We consider linear systems on toric varieties of any dimension, 
with invariant base points, giving a characterization of special linear systems. 
We then make a new conjecture for linear systems on rational surfaces.
\end{abstract}
\maketitle

\section*{Introduction} 
Let us take the projective plane $\p^2$ and 
let us consider the linear system of curves of degree $d$ having some points 
of fixed multiplicity. The expected dimension of such systems is 
the dimension of the space of degree $d$ polynomials minus the conditions
imposed by the multiple points. 
The systems whose dimension is bigger than the expected one are called
{\em special systems}.

There exists a conjecture due to Hirschowitz (see \cite{hi}), characterizing special 
linear systems on $\p^2$, which has been proved in some special case 
\cite{cm,cm2,mig1,ev}
There exists also an analogous conjecture for rational ruled surfaces \cite{la}.

Concerning the dimension of linear systems on varieties of bigger dimension, with 
some points of fixed multiplicity, very few is known. 

In this article we study this problem in the case of toric varieties, considering
only points which are invariant under the action of the maximal torus. 
In particular we take a smooth $n$-dimensional toric variety
$X$, an ample divisor $D$ on it, $r$ equivariant points $p_1,p_2,\ldots ,p_r$ 
and $r$ non negative integers $m_1,m_2,\ldots ,m_r$. We denote by 
$\ls(D,m_1,\ldots ,m_r)$
the linear system of the divisors linearly equivalent to $D$ and passing through $p_i$ with 
fixed multiplicity $m_i$. 

Then we prove the following:  
\begin{mainthm}
The linear system $\ls(D,m_1,\ldots ,m_r)$ is special if and only if there exists an 
equivariant curve 
$C$ passing through two of the $p_i$'s, such that $\ls(D,m_1,\ldots ,m_r)\cdot C\leq -2$.
\end{mainthm}

The paper is organized as follows: in Section $1$ we recall some definitions and we fix
some notations (see for instance \cite{oda} for a complete reference), while 
in Section $2$ we state and prove
our main theorem and we formulate a new conjecture for linear systems on rational surfaces.

\section{Definitions and Notations}

In what follows we will adopt the notation of \cite{oda}, Chapter 2.

Let $N\cong\z^n$ and let $M$ be its dual
$\z$-module; we denote by $N_{\rea}$ (resp. $M_{\rea}$) the $\rea$-module $N\otimes\rea$
(resp. $M\otimes\rea$). Given a fan $\Delta$ in $N$, the corresponding
$n$-dimensional toric variety $X$ is denoted by $T_N emb(\Delta)$. 
In this paper we consider only smooth toric varieties defined by nonsingular, complete
fans. We denote by $\sigma_1,\ldots,\sigma_k$ the $n$-dimensional
cones of $\Delta$ corresponding to the $T_N$-invariant points of $X$, $p_1,\ldots,p_k$.

Given an ample divisor $D$ on $X$, let $h\in SF(N,\Delta)$ be its $\Delta$-linear support
function. The convex polytope in $M_{\rea}$ representing $H^0(X,\oc_X(D))$ is denoted by
$\square_h=\{m\in M_{\rea}\mid \langle m,n\rangle\geq h(n),\forall n
\in N_{\rea}\}$.

\vskip .2truecm

Let $V_i$ be the vertex of $\square_h$ corresponding to the element of $|D|$ that does not 
pass through the $T_N$-invariant point $p_i$, for $i=1,\ldots, r$. We denote by 
$C_{j,k}$ the $T_N$-invariant curve joining the points $p_j$ and $p_k$ (if it exists), by
$\e_{j,k}$ the edge of $\square_h$ corresponding to the elements of $|D|$ that does not contain
$C_{j,k}$ (i.e. the edge joining $V_j$ and $V_k$), and by $N_{j,k}$ the number of integer
points lying on $e_{j,k}$.

In order to explain the notations, we are going to see the example of
the system $|3H+2F|$ on the Hirzebruch surface $\f_3$ (we denote by $H$
a rational $3$-curve and by $F$ the fiber).

\vspace{.3cm}
{
\unitlength=.7pt
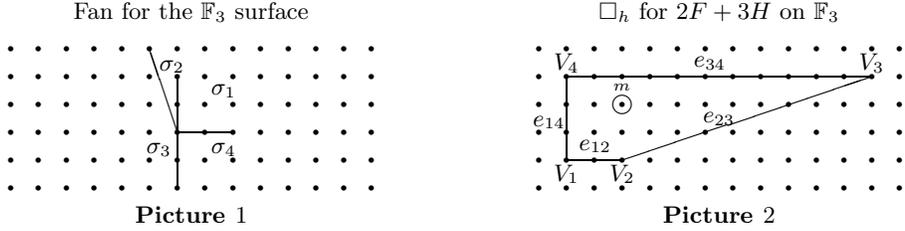
\begin{picture}(210,135)(-2,0)
\put(0,0){\makebox(195,25)[t]{\small {\bf Picture $1$}}}
\put(0,110){\makebox(195,25)[t]{\small Fan for the $\f_3$ surface}}
\multiput(0,35)(15,0){14}{\circle*{3}}
\multiput(0,50)(15,0){14}{\circle*{3}}
\multiput(0,65)(15,0){14}{\circle*{3}}
\multiput(0,80)(15,0){14}{\circle*{3}}
\multiput(0,95)(15,0){14}{\circle*{3}}
\multiput(0,110)(15,0){14}{\circle*{3}}
\put(90,65){\line(1,0){30}}
\put(90,65){\line(0,1){30}}
\put(90,65){\line(0,-1){30}}
\put(90,65){\line(-1,3){15}}
\put(100,80){\makebox(30,15){\small$\sigma_1$}}
\put(72,93){\makebox(30,15){\small$\sigma_2$}}
\put(65,48){\makebox(30,15){\small$\sigma_3$}}
\put(100,48){\makebox(30,15){\small$\sigma_4$}}
\end{picture}}
\hspace{1.5cm}
{
\unitlength=.7pt
\begin{picture}(210,135)(-2,0)
\put(0,0){\makebox(195,25)[t]{\small {\bf Picture $2$}}}
\put(0,110){\makebox(195,25)[t]{\small $\square_h$ for $2F+3H$ on $\f_3$}}
\put(0,35){\makebox(30,15){\small$V_1$}}
\put(30,35){\makebox(30,15){\small$V_2$}}
\put(165,95){\makebox(30,15){\small$V_3$}}
\put(0,95){\makebox(30,15){\small$V_4$}}

\put(23,50){\makebox(15,15){\small $e_{12}$}}
\put(90,65){\makebox(15,15){\small $e_{23}$}}
\put(85,95){\makebox(15,15){\small $e_{34}$}}
\put(-2,63){\makebox(15,15){\small $e_{14}$}}

\multiput(0,35)(15,0){14}{\circle*{3}}
\multiput(0,50)(15,0){14}{\circle*{3}}
\multiput(0,65)(15,0){14}{\circle*{3}}
\multiput(0,80)(15,0){14}{\circle*{3}}
\multiput(0,95)(15,0){14}{\circle*{3}}
\multiput(0,110)(15,0){14}{\circle*{3}}
\put(15,50){\line(1,0){30}}
\put(45,50){\line(3,1){135}}
\put(15,95){\line(1,0){165}}
\put(15,50){\line(0,1){45}}
\put(45,80){\circle{10}}
\put(40,82){\makebox(10,15){\tiny $m$}}
\end{picture}} \\

The first picture represents the complete fan of a $\f_3$ surface.
In the second picture the integer points of $\square_h$ represent the elements 
of $H^0(2F+3H)$. Here for instance we have that $e_{1,4}$ corresponds a fiber 
$F$, and $N_{1,4}=4$.

To each $m\in\square_h$ there is associated 
an element $s_m\in H^0(2F+3H)$. 
To calculate the multiplicity of $s_m$ at $p_i$ we translate the polytope $\square_h$ 
in such a way that $V_i$ coincides with the vertex of the dual cone 
$\sigma_i^{\vee}$, and hence 
we express $m$ as a combination of the generators of $\sigma_i^{\vee}$. The multiplicity of 
$s_m$ at $p_i$ is just the sum of the coefficients in that combination.
For instance, if $m$ is the point in the picture we have
that the multiplicities of $s_m$ at $p_1,\ldots,p_4$ are $4,8,7,3$ respectively.

\vskip .2truecm

Let us fix the $T_N$-invariant point $p_i$ and hence the vertex $V_i$. There exist $n$ edges 
$e_{i,k_1},\ldots,e_{i,k_n}$, passing through $V_i$ (because the variety $X$ is smooth). 
Let us denote by $H_i$ the hyperplane meeting the edges $e_{i,k_j}$ ($j=1,\ldots,n$) in 
the $m_i$th integer point, starting from $V_i$ (we consider the extension of the edge
if the number $N_{i,k_j}$ of integer points on it is smaller than $m_i$).
To impose that the system passes through the point $p_i$ with multiplicity $m_i$ is equivalent
to cut from $\square_h$ its intersection with the half-space determined by $H_i$ and 
containing the vertex $V_i$. 
It is clear that the number of integer points we are cutting is exactly $\binom{m_i+n-1}{n}$
(which are the conditions imposed by a point of multiplicity $m_i$) if and only if the number
of integer points lying on each edge $e_{i,k_j}$, 
is at least $m_i$ (i.e. $N_{i,k_j}\geq m_i$), for $j=1,\ldots,n$.

\vskip .3truecm

{
\unitlength=.7pt
\begin{picture}(210,135)(-2,-35)
\put(0,-35){\makebox(195,25)[t]{\small {\bf Picture $3$}}}
\put(0,75){\makebox(195,25)[t]{\small $\mid 2F+3H-5p_4\mid$ on $\f_3$}}
\put(20,45){\makebox(30,15){\tiny $m_4=5$}}
\put(0,0){\makebox(30,15){\small$V_1$}}
\put(30,0){\makebox(30,15){\small$V_2$}}
\put(165,60){\makebox(30,15){\small$V_3$}}
\put(0,60){\makebox(30,15){\small$V_4$}}

\multiput(0,0)(15,0){14}{\circle*{3}}

\put(0,15){\circle*{3}}
\multiput(15,15)(15,0){2}{\circle{4}}
\multiput(45,15)(15,0){11}{\circle*{3}}

\put(0,30){\circle*{3}}
\multiput(15,30)(15,0){3}{\circle{4}}
\multiput(60,30)(15,0){10}{\circle*{3}}

\put(0,45){\circle*{3}}
\multiput(15,45)(15,0){4}{\circle{4}}
\multiput(75,45)(15,0){4}{\circle*{3}}
\put(135,45){\circle*{3}}  
\multiput(150,45)(15,0){4}{\circle*{3}}

\put(0,60){\circle*{3}}
\multiput(15,60)(15,0){5}{\circle{4}}
\multiput(90,60)(15,0){5}{\circle*{3}}
\multiput(165,60)(15,0){2}{\circle*{3}} 
\put(195,60){\circle*{3}}

\multiput(0,75)(15,0){14}{\circle*{3}}

\put(15,15){\line(1,0){30}}
\put(45,15){\line(3,1){135}}
\put(15,60){\line(1,0){165}}
\put(15,15){\line(0,1){45}}
\put(20,5){\line(1,1){65}}
\end{picture}}
\hspace{1.5cm}
{
\unitlength=.7pt
\begin{picture}(210,135)(-2,-35)
\put(0,-35){\makebox(195,25)[t]{\small {\bf Picture $4$}}}
\put(0,75){\makebox(195,25)[t]{\small $\mid 2F+3H-2p_1-3p_4\mid$ on $\f_3$}}
\put(20,45){\makebox(30,15){\tiny $m_4=3$}}
\put(20,15){\makebox(30,15){\tiny $m_1=2$}}
\put(0,0){\makebox(30,15){\small$V_1$}}
\put(30,0){\makebox(30,15){\small$V_2$}}
\put(165,60){\makebox(30,15){\small$V_3$}}
\put(0,60){\makebox(30,15){\small$V_4$}}
\multiput(0,0)(15,0){14}{\circle*{3}}

\put(0,15){\circle*{3}}
\multiput(15,15)(15,0){2}{\circle{4}}
\multiput(45,15)(15,0){11}{\circle*{3}}

\put(0,30){\circle*{3}}
\put(15,30){\circle{4}}
\multiput(30,30)(15,0){12}{\circle*{3}}

\put(0,45){\circle*{3}}
\multiput(15,45)(15,0){2}{\circle{4}}
\multiput(45,45)(15,0){11}{\circle*{3}}

\put(0,60){\circle*{3}}
\multiput(15,60)(15,0){3}{\circle{4}}
\multiput(60,60)(15,0){10}{\circle*{3}}

\multiput(0,75)(15,0){14}{\circle*{3}}
\put(15,15){\line(1,0){30}}
\put(45,15){\line(3,1){135}}
\put(15,60){\line(1,0){165}}
\put(15,15){\line(0,1){45}}
\put(5,20){\line(1,1){50}}
\put(5,40){\line(1,-1){35}}
\end{picture}}
\vspace{.5cm}

In Picture $3$ and $4$ we represent the special systems 
$\mid 2F+3H-5p_4\mid$ and $\mid 2F+3H-2p_1-3p_4\mid$ on $\f_3$. 
The speciality of the first system is due to the conditions imposed on $p_4$,
while the speciality of the second one is due to the overlapping of conditions imposed 
on $p_1$ and $p_4$.

\vskip .2truecm

\begin{notation}
Throughout the paper, by abuse of notation, we will use the same symbol 
to designate a linear system and the corresponding sheaf. 
\end{notation}

\begin{defi}
We say that a linear system $\ls$ on a smooth variety $X$ is {\em special} if $h^1(X,\ls)\neq 0$.
\end{defi}

Given an ample, non-special divisor $D$ on $X$,
(i.e. such that $\oc_X(D)$ is non-special), $r$ points 
$p_1,\ldots,p_r\in X$, and $r$ non negative integers $m_1,\ldots,m_r$,
we denote by $\ls(D,m_1,\ldots ,m_r)$ the linear system of divisors linearly equivalent 
to $D$ and passing through $p_i$ with multiplicity $m_i$. Therefore, if we 
denote by $Z$ the $0$-dimensional scheme of multiple points, the linear system corresponds
to the sheaf $\oc_X(D)\otimes\I_Z$.

\begin{defi}
The {\em virtual dimension} of the system $\ls(D,m_1,\ldots ,m_r)$, is denoted by
$v(\ls(D,m_1,\ldots ,m_r))$, and is defined to be the difference between the 
dimension of the complete system $|D|$ and the conditions imposed by the multiple points
i.e.
\begin{equation*}
v(\ls(D,m_1,\ldots ,m_r)) = h^0(X,\oc_X(D))-h^0(Z,\oc_Z(D))-1.
\end{equation*}
\end{defi}

From the exact sequence
\begin{equation*}
0\rightarrow \I_Z\rightarrow \oc_X\rightarrow\oc_Z\rightarrow 0,
\end{equation*}
tensoring with $\oc_X(D)$ and taking cohomology, we obtain
\begin{equation*}
0\rightarrow \h^0(X,\ls)\rightarrow \h^0(X,\oc_X(D))\rightarrow
\h^0(Z,\oc_Z\otimes\oc_X(D)) \rightarrow \h^1(X,\ls) \rightarrow 0,
\end{equation*}
where, for simplicity of notation, we put $\ls:=\ls(D,m_1,\ldots ,m_r)$.
Therefore
\begin{equation*}
v(\ls) = h^0(X,\ls)-h^1(X,\ls)-1,
\end{equation*}
and, since the effective dimension of $\ls$ is $h^0(X,\ls)-1$, the system is special 
if and only if its virtual dimension is smaller than the effective one.

\section{Main Theorem and Conjecture}

Let $X$ be a smooth $n$-dimensional toric variety, and let $D$ be an ample divisor
on it. We fix $r$ points $p_1,p_2,\ldots ,p_r$ on $X$, $T_N$-invariant,
$r$ non negative integers $m_1,m_2,\ldots ,m_r$, and we consider the linear system 
$\ls(D,m_1,\ldots ,m_r)$.

Let $C$ be a curve on $X$, passing through some of the points $p_i$, say $p_1,\ldots, p_s$, 
$s\leq r$. 
We define the intersection of the system $\ls(D,m_1,\ldots,m_r)$ and the curve $C$ as the 
intersection product of their strict transforms $\ls'$ and $C'$ on $X'$, the blow up of 
$X$ along the points $p_1,\ldots, p_s$. 
Therefore, if we denote by $E_i\cong\p^{n-1}$ the exceptional divisor corresponding 
to the point $p_i$, and by $C'=\overline{\pi^{-1}(C\setminus\{p_1\ldots p_r\})}$, 
since $C'\cdot E_i=1$ for each $i=1\ldots r$, then we have the following formula:

\begin{equation}\label{blow}
\begin{array}{ll}
C\cdot\ls(D,m_1,\ldots ,m_r) & = C'\cdot\ls'\\
& = C'\cdot(\pi^*D-\sum_{i=1}^r m_iE_i)\\
& = C\cdot D-\sum_{i=1}^s m_i.
\end{array}
\end{equation}

In order to simplify the proof of our main theorem we are going to state and prove
two lemmas.

\begin{lemma}\label{lem1}
The system $\ls(D,m_1,\ldots ,m_r)$ is special if and only if there exist two points, say
$p_1$ and $p_2$, such that the system $\ls(D,m_1,m_2)$ is special.
\end{lemma}

\begin{proof}
The system $\ls(D,m_1,\ldots,m_r)$ is special if and only if one of
the following hold.

On one hand, there can exist a point, say  $p_1$, such that $m_1$ is bigger than $N_{1,k_j}$
for some index $j\in\{1,\ldots,n\}$, which means that the system $\ls(D,m_1)$ (and hence
also $\ls(D,m_1,m_2)$) is already special.

On the other hand, we have speciality if some of the half-spaces we are cutting
have a nonempty intersection inside the polytope $\square_h$.
If this is the case, in particular there must exist 
two of the half-spaces, corresponding to say $H_1$ and $H_2$, intersecting inside $\square_h$, 
which is equivalent to say that the system $\ls(D,m_1,m_2)$ is already special. 
\end{proof}

\begin{lemma}\label{lem2}
Let $D$ and $C_{j,k}$ be as above. Then the following equality holds:
\begin{equation}
D\cdot C_{j,k}=N_{j,k}-1.
\end{equation}
\end{lemma}

\begin{proof}
Let us write $D=\sum_{i=1}^m\alpha_i D_i$, where $D_i$ is the $T_N$-invariant divisor on $X$,
corresponding to the $1$-cone $v_i$, and let us take the curve $C_{j,k}$, 
corresponding to the $(n-1)$-cone $\tau=\langle v_1,\ldots,v_{n-1}\rangle$. 
Let us denote by $v_n$ and $v_{n+1}$ the $1$-cones which complete $\tau$
to the two $n$-cones containing it. 
Therefore the intersection product of $C_{j,k}$ with $D$ depends only on its intersection with 
the $D_i$, for $i=1,\ldots,n+1$ (see \cite{oda}, Chapter $2$). Actually we have 
\begin{equation*}
D\cdot C_{j,k}=\sum_{i=1}^{n+1}\alpha_i D_i\cdot C_{j,k}.
\end{equation*} 
Since $\langle v_1,\ldots,v_{n-1},v_n\rangle$ and $\langle v_1,\ldots,v_{n-1},v_{n+1}\rangle$
are two $n$-cones, it follows that $D_n\cdot C_{j,k}=D_{n+1}\cdot C_{j,k}=1$. Let us then 
calculate $D_i\cdot C_{j,k}$, for $i=1,\ldots, n-1$.

Let us remark that, via unimodular transformation, we can always suppose that $v_i=e_i$, for 
$i=1,\ldots,n$, (i.e. the canonical base for $\rea^n$. Since  
$|v_1,\- \ldots,\- v_{n-1},\- v_{n+1}|=-|v_1,\ldots,v_{n-1},v_n|=-1$, the last coordinate 
of $v_{n+1}$ must be $-1$ and hence we can write $v_{n+1}=(-\gamma_1,\ldots,-\gamma_{n-1},-1)$. 
Therefore the following equality holds
\begin{equation*}
v_n+v_{n+1}+\gamma_1v_1+\ldots+\gamma_{n-1}v_{n-1}=0,
\end{equation*}
which gives (see \cite{oda}) $D_i\cdot C_{j,k}=\gamma_i$, for $i=1,\ldots, n-1$, and hence
\begin{equation}\label{inters}
D\cdot C_{j,k}=\sum_{i=1}^{n-1}\alpha_i\gamma_i+\alpha_n+\alpha_{n+1}.
\end{equation}
Let us computate now $N_{j,k}$, i.e. the number of integer points on the edge $e_{j,k}$ of 
$\square_h$. This edge lies on the line $r=\{x_i=-\alpha_i,\ i=1,\ldots n-1\}$. In particular
it is parallel to the $x_n$-axis and hence $N_{j,k}=l+1$, where $l$ is the lenght of $e_{j,k}$.
Therefore, since the edge is cut by the two hyperplanes 
\begin{equation*}
\begin{array}{ll}
H & = \{w\mid \langle w,v_n\rangle=-\alpha_n\}\\
  & = \{x_n=-\alpha_n\}
\end{array}
\end{equation*}
and
\begin{equation*}
\begin{array}{ll}
H' & = \{w'\mid \langle w',v_{n+1}\rangle=-\alpha_{n+1}\}\\
   & = \{\gamma_1x_1+\ldots\gamma_{n-1}x_{n-1}+x_n=\alpha_{n+1}\},
\end{array}
\end{equation*}
its endpoints are $p=(-\alpha_1,\ \ldots,\ -\alpha_{n-1},\ -\alpha_n)$ 
and $p'=(-\alpha_1,\ \ldots,\ -\alpha_{n-1},$\  
$\sum_{i=1}^{n-1}\ \gamma_i\alpha_i+\alpha_{n+1})$. Therefore, recalling
\eqref{inters},
\begin{equation*}
\begin{array}{ll}
N_{j,k}-1 & = l\\
          & = \sum_{i=1}^{n-1}\gamma_i\alpha_i+\alpha_{n+1}+\alpha_n\\
          & = D\cdot C_{j,k}.
\end{array}
\end{equation*}
\end{proof}

\begin{mainthm}
The linear system $\ls(D,m_1,\ldots ,m_r)$ is special if and only if there exists a 
$T_N$-invariant curve 
$C$ passing through two of the $p_i$'s, such that $\ls(D,m_1,\ldots ,m_r)\cdot C\leq -2$.
\end{mainthm}
\begin{proof}

Because of Lemma \ref{lem1} we can fix our attention on the system $\ls(D,m_1,m_2)$, 
relating its speciality with the number $N_{1,2}$. 

In fact, if the speciality of our system is due to only one point, say $p_1$, 
then there exists an edge, say $e_{1,2}$ such that $m_1\geq N_{1,2}+1$,
as we have seen in the proof of previous Lemma \ref{lem1}.
On the other hand, if the speciality is due to both $p_1$ and $p_2$, then
the two half-spaces determined by $H_1$ and $H_2$ intersect inside $\square_h$.
In this case their intersection contains at least one point of the edge $e_{1,2}$. 

Therefore, in both cases we have that the system is special if and only if the 
following equality holds:
\begin{equation}\label{m1m2}
m_1+m_2\geq N_{1,2}+1.
\end{equation}

From the intersection formula 
\eqref{blow} and from Lemma \ref{lem2},
\begin{equation*}
\begin{array}{ll}
\ls(D,m_1,m_2)\cdot C_{1,2} & = D\cdot C_{1,2}-(m_1+m_2)\\
                            & = N_{1,2}-1-(m_1+m_2).\\
\end{array}
\end{equation*}
and hence, from \eqref{m1m2}, the system is special if and only if $\ls(D,m_1,m_2)\cdot C_{1,2}
\leq -2$.
\end{proof}

\vskip .2truecm

\begin{example}
Let us consider the Hirzebruch surface $\f_3$. This is the toric surface 
associated to the complete two-dimensional fan of Picture $1$.
There exist $4$ invariant 1-dimensional varieties, namely the $(-3)$-curve $E$, 
one rational $3$-curve $H$ and two fibers $F_1$ and $F_2$. 
Their intersections give rise to $4$ invariant points, which we denote
by $p_i,\ i=1,\ldots,4$ as in Picture $2$.

Let us consider now the linear system $|3H+2F-5p_4|$. Blowing up along $p_4$ one 
obtains the system $|3\pi^*H+2\pi^*F-5E_4|$, where $E_4$ is the exceptional divisor of
the blowing up $\pi$. Now observe that $|\pi^*F-E_4|$ is a curve of genus $0$ and 
that $(3\pi^*H+2\pi^*F-5E_4)\cdot (\pi^*F-E_4)=-2$, hence by the remark below, the system
must be special (as we can see geometrically from Picture $3$). 
In a similar way also the system $|3H+2F-2p_1-3p_4|$ 
is special since, blowing up $p_1$ and $p_4$ one has 
$(3\pi^*H+2\pi^*F-2E_1-3E_4)\cdot (\pi^*F-E_1-E_4)=-2$ where $|\pi^*F-E_1-E_4|$ is the fiber 
through $p_1,p_4$ (see Picture $4$).
\end{example}

\begin{remark}
Observe that if $S$ is a rational surface, $D$ is an ample, non special divisor, 
and $\ls$ is a linear system of curves linearly equivalent to $D$ and
such that there exists a rational, irreducible  curve 
$C\subset S$ with $\ls\cdot C\leq -2$, then $\ls$ is special.

In fact, the effective dimension of $\ls$ equals that of $\ls - C$, and we are going to see
that the virtual dimension $v(\ls)$ is smaller than $v(\ls-C)$
(and hence, in particular, $v(\ls)$ is smaller than the effective dimension of $\ls$,
and the system is special). 
We recall that $v(\ls)=h^0(\ls)-h^1(\ls)-1$ and,
being $S$ rational, $h^2(\ls)=0$, which implies $v(\ls)=\chi(\ls)-1$.
By Riemann-Roch, $\chi(\ls)=(\ls^2-\ls\cdot K_S)/2+1$. 
Therefore $v(\ls)=(\ls^2-\ls\cdot K_S)/2$ and  
$v(\ls-C)=((\ls-C)^2-(\ls-C)\cdot K_S)/2 = v(\ls)+g(C)-1-\ls\cdot C$.
The rationality of $C$ implies that
\[
v(\ls)=v(\ls-C)+\ls\cdot C+1.
\]
Hence, if $\ls\cdot C\leq -2$, we have that $v(\ls) < v(\ls-C)$.
\end{remark}

Consider the Hirzebruch surface $\f_6$ and the linear system $\ls(\oc_{\f_6}(0,4),3^{11})$. 
The virtual dimension of this system is $-2$.
Let $C$ be an element of $\ls(\oc_{\f_6}(2,1),1^{11})$, then $C$ is rational and 
$\ls\cdot C=-1$. Let $\Gamma_6\in \ls(\oc_{\f_6}(-6,1))$ be the $(-6)$-curve. 
Then we have the following:

\[
\begin{array}{ccl}
\ls \cdot C			& = & -1 \\
(\ls - C) \cdot \Gamma_6	& = & -2 \\
(\ls - C - \Gamma_6) \cdot C	& = & -2 \\
\end{array}
\]

Now, since 
$3C + \Gamma_6 = 3\ls(\oc_{\f_6}(2,1),1^{11}) + \ls(\oc_{\f_6}(-6,1)) = 
\ls(\oc_{\f_6}(0,4),3^{11})$, then the initial system is not empty and therefore special.
Observe that there exists no curve $E$ on $\f_6$ such that $\ls \cdot E\leq -2$.\\

This motivates the following procedure in the case of a smooth rational surface $X$.\\

{\bf Step 1} - Let $\ls$ be a non-empty linear system on a rational surface $X$. 
If there exists a rational curve $C$ on $X$ such that $\ls\cdot C\leq -1$, replace $\ls$ with
$\ls - C$ and restart from Step 1. \\
{\bf Step 2} - If the last system has virtual dimension bigger than that of the initial one 
then we say that $\ls$ is a {\em $(-1)$-special system}. \\

Considering the remark above and our main theorem, we conjecture that:

\begin{con}\label{surf}
Let $D$ be an ample, non-special divisor on a rational surface $S$, and let $p_1,\ldots,p_r$ be
$r$ points on $S$. Then the linear system
$\ls=\ls(D,m_1,\ldots,m_r)$ is special if and only if it is $(-1)$-special.

\end{con}

\bibliographystyle{amsplain}

\begin{thebibliography}{1}

\bibitem{cm}
Ciro Ciliberto and Rick Miranda, \emph{Degenerations of planar linear systems},
  J. Reine Angew. Math. \textbf{501} (1998), 191--220.

\bibitem{cm2}
\bysame, \emph{{Linear systems of plane curves with base points of equal
  multiplicity.}}, Trans. Am. Math. Soc. \textbf{352} (2000), no.~9, 4037--4050
  (English).

\bibitem{hi}
Andr{\'e} Hirschowitz, \emph{Une conjecture pour la cohomologie des diviseurs
  sur les surfaces rationnelles g\'en\'eriques}, J. Reine Angew. Math.
  \textbf{397} (1989), 208--213.

\bibitem{la}
Antonio Laface, \emph{On linear systems of curves on rational scrolls}, Geom.
  Dedicata \textbf{90} (2002), no.~1, 127--144.

\bibitem{ev}
Evain Laurent, \emph{La fonction de {H}ilbert de la r\'eunion de $4\sp h$ gros
  points g\'en\'eriques de ${\bf {p}}\sp 2$ de m\^eme multiplicit\'e}, J.
  Algebraic Geom. \textbf{8} (1999), no.~4, 787--796.

\bibitem{mig1}
Thierry Mignon, \emph{Syst\`emes de courbes planes \`a singularit\'es
  impos\'ees: le cas des multiplicit\'es inf\'erieures ou \'egales \`a quatre},
  J. Pure Appl. Algebra \textbf{151} (2000), no.~2, 173--195.

\bibitem{oda}
Tadao Oda, \emph{{Convex bodies and algebraic geometry. An introduction to the
  theory of toric varieties.}}, {Ergebnisse der Mathematik und ihrer
  Grenzgebiete. 3. Folge, Bd. 15. Berlin etc.: Springer-Verlag. VIII, 212 p.;
  DM 148.00 }, 1988 (English).

\end{thebibliography}

\providecommand{\bysame}{\leavevmode\hbox to3em{\hrulefill}\thinspace}

\end{document}